\documentclass[12pt]{amsart}

\usepackage{amssymb}

\usepackage{enumerate}

\usepackage{graphicx}

\makeatletter
\@namedef{subjclassname@2010}{%
  \textup{2010} Mathematics Subject Classification}
\makeatother

 \usepackage[T1]{fontenc}

\newtheorem{thm}{Theorem}[section]

\newtheorem{lem}[thm]{Lemma}

\theoremstyle{definition}

\newtheorem*{rqes}{Remarks}

\numberwithin{equation}{section}

\def\m{\mathcal}

\def\C{\mathbb{C}}
\def\c2{\mathbb{C}^2}
\def\R{\mathbb{R}}
\def\N{\mathbb{N}}

\def\N{\mathbb{N}}

\def\1{\bold{1}}

\def\f{\varphi}

\newcommand \W {\Omega}

\newcommand \mE {\mathcal E}

\newcommand \Sub {\Subset}

\begin{document}


\baselineskip=17pt



\title[Weak Solutions To Complex Monge-Amp\`ere  Equation]
{Weak Solutions To Complex Monge-Amp\`ere  Equation on hyperconvex domains}

\author[S. Benelkourchi]{Slimane Benelkourchi}
\address{Universit\'e de Montr\'{e}al,
Pavillon 3744, rue Jean-Brillant,
Montr\'{e}al QC  H3C 3J7, Canada.}
\email{slimane.benelkourchi@umontreal.ca}
\date{}
\begin{abstract}
 We show a very general existence theorem to the complex
Monge-Amp\`ere type equation on hyperconvex domains.
\end{abstract}



\maketitle

\section{Introduction}
 Let $\W$
 be a bounded hyperconvex domain in $\C ^n$ and $F$  a nonnegative function defined on
$\R \times \W.$
In the present note, we shall consider the existence and uniqueness
 of weak solution of
  the complex Monge-Amp\`ere type equation
\begin{equation}\label{local}
(dd^c u)^n = F(u , \cdot) d \mu
\end{equation}
 where $u$ is plurisubharmonic on
  $\W $ and  $\mu$ is a nonnegative measure.
 This problem has been studied extensively by various authors,
  see for example \cite{BT0}, \cite{BT1},
 \cite{Ce8}, \cite{Ce4}, \cite{Ce98} \cite{CK6}, \cite{CK94}, \cite{CY}, \cite{K98}, \cite{K05}, \cite{K0}
, \cite{Y}... and reference therein for further information about Complex Monge-Amp\`ere
equations.

It was first considered by Bedford and Taylor in \cite{BT79}. In connection with the problem of finding
complete K\"ahler-Einstein metrics on pseudoconvex domains, Cheng and Yau \cite{CY} treated the case
$F(t, z) = e^{Kt} f(z)$. More recently, Czy\.z \cite{Cz9} treated the case $F$ bounded by a function
independent  of the first variable and $\mu$ is the Monge-Amp\`ere of a plurisubharmonic
 function $v,$ generalizing some results of Cegrell \cite{Ce84}, Ko\l odziej \cite{K0} and Cegrell and
Ko\l odziej \cite{CK6} \cite{CK94}. In this paper we will consider a  more general case.
With notations introduced in the next section, our main result  is stated as follows.\\

\noindent{\it {\bf Main Theorem }
Let $\W$ be a bounded hyperconex domain and
$\mu $ be a nonnegative measure which vanishes on all pluripolar subsets of $\W. $
Assume that $F : \R \times \W \to [0, +\infty )$ is a measurable function such that:

1) For all $z\in \W $ the function $t \mapsto F(t, z)$ is continuous and  nondecreasing;

2) For all $t\in \R $, the function $z\mapsto F(t,z) $ belongs to $L^1_{loc}(\W, \mu)$;

3) There exists a function $v_0 \in \m N^a $ which is a subsolution to (\ref{local}) i.e
$$
(dd^c v_0)^n \ge F(v_0 , \cdot) d \mu .
$$

Then for  any maximal function $f \in \m E$ there exists a uniquely determined function
$u\in \m N^a(f)$ solution to the complex Monge-Amp\`ere equation
\begin{equation*}
(dd^c u)^n = F(u , \cdot) d \mu .
\end{equation*}
}

 Note that the solution, as we will see in the proof,
  is given by the following upper envelope of all subsolutions;
 $$
 u = \sup \left\lbrace v \in \m E(\W) ; v \le f  \ \text{and}\
  (dd^c v)^n \ge F(v , \cdot) d \mu \right\rbrace,
 $$
 where $\m E(\W) $ is the  set of non-positive plurisubharmonic functions defined on
$\W$
 for which the
complex Monge-Amp\`ere operator is well defined as nonnegative measure
(a precise definition will be given shortly).
 \section{Background and Definitions}
 Recall that $\W \Sub  \C^n ,\  n \ge 1$ is a bounded hyperconvex domain if it is a bounded,
connected, and open set such that there exists a bounded plurisubharmonic function
$\rho : \W \to (-\infty , 0)$ such that the closure of the set
 $\left\lbrace z\in \W  :\  \rho(z) < c\right\rbrace$ is compact in
$\W,$ for every $ c \in  (-\infty, 0).$ We denote by $PSH(\W)$
 the family of plurisubharmonic
functions defined on $\W.$

We say that a bounded plurisubharmonic function $\f $ defined on $\W$
 belongs to $\m E_0 $
if $ lim _{z \to \zeta}  \f(z) = 0,$ for every $\zeta \in \partial\W,$ and
$\int_\W  (dd^c \f)^n < +\infty.$ See \cite{Ce4} for details.

 Let the class $\m E (\W)$
 be the set of plurisubharmonic functions
 $u $
such that for all $z_0 \in \Omega $,
 there exists a neighborhood $V_{z_0}$ of
$z_0$ and $u_j \in {\mathcal E}_0(\Omega)$ a decreasing sequence
which converges towards $u$ in $V_{z_0}$ and satisfies
$$\sup_j
\int_{\Omega} (dd^c u_j)^n <+\infty.$$
 U.Cegrell
\cite{Ce4} has shown  that the operator $(dd^c \cdot )^n$ is well defined on
$\mE(\Omega),$ continuous under decreasing limits and the class
$\mE(\Omega)$ is stable under taking maximum i.e. if $u\in \mE(\W) $
 and $v\in PSH^-(\W)$ then $\max (u , v) \in \mE(\W).$ $\m E(\W)$ is the largest
class with these properties (Theorem 4.5 in \cite{Ce4}). The class
$\mE(\Omega)$ has been further characterized by Z.B\l{}ocki \cite{Bl1}, \cite{Bl2}.

The class ${\mathcal F}(\Omega)$ is the ``global version'' of
 $\mE(\Omega)$:
a function $u$ belongs to ${\mathcal F}(\Omega)$ iff there exists
 a decreasing sequence
 $u_j \in {\mathcal E}_0(\Omega)$
converging towards $u$ {\it in all of } $\Omega$, which satisfies
$\sup_j \int_{\Omega} (dd^c u_j)^n<+\infty.$  Furthermore characterizations are given in \cite{B1} \cite{BGZ2}.


 Define $\m N(\W)$ the family of all function $u \in \m E(\W) $ which satisfies:
 if $v \in PSH(\W) $ is  maximal and $u\le v$ then $v\ge 0,$
i.e. the smallest
maximal psh function above $u$ is null.
In fact, this class is the analogous of potentials
for subharmonic functions (see \cite{Ce8} for more details).

The class $\m F^ a(\W)$ (resp. $\m N^a(\W)$, $\m E^a (\W)$
 $\cdots$)
 is the set of functions $u \in \m F (\W)$
 (resp. $u \in \m N (\W)$,
 $u \in \m E (\W)$ $\cdots$)
  whose Monge-
Amp\`ere measure $(dd^cu)^n$ is absolutely continuous with respect to
capacity i.e. it does not charge pluripolar sets.

Finally, for $f\in \m E,$ we denote by $\m N (f) $ (resp. $\m F (f)$)
the family of those  $u\in PSH(\W)$  such that there exists a
function
$\f \in \m N $ (resp. $\f \in \m F$)  satisfying the following inequality
$$
\f (z) + f(z) \le u (z) \le f(z) \quad \forall z \in \W.
$$

We shall use  repeatedly the following
 well known comparison principle
from \cite{BT1}
as well as its generalizations to the class $\m N(f)$ (cf  \cite{ACCP} \cite{Ce8}).
\begin{thm}[ \cite{ACCP} \cite{BT1} \cite{Ce8}]
Let $ f \in \m E(\W)$ be a maximal function and $u, \ v,\ \in \m N(f)$  be such that $(dd^c u)^n $
vanishes on all pluripolar sets in $\W.$ 
 Then
$$
\int _{(u<v)} (dd^c v)^n \le \int _{(u<v)} (dd^c u)^n .
$$
Furthermore if $ (dd^c u)^n=  (dd^c v)^n$ then $u=v.$
\end{thm}

\section{Proof of Main Theorem}
\begin{lem}[Stability] \label{stability}Let $\mu $
 be a finite nonnegative measure which vanishes on all pluripolar subsets of $\W$
and $f\in \m E (\W) $ be a maximal function. Fix a function
 $v_0 \in \m E(\W).$ Then for any $u_j ,\ u \in \m N^a(f) $
solutions to
$$
(dd^c u_j)^n =h_j d\mu , \quad (dd^c u)^n =h d\mu
$$
such that $0\le h d\mu , h_j d\mu\le (dd^c v_0)^n $
 and $h_j d\mu \to hd\mu $ as measures, we have that $u_j $ converges towards $u$ weakly.
\end{lem}
The statement of the lemma fails if no control on
 the complex Monge-Amp\`ere measures is assumed (see \cite{CK94}).
\begin{proof}
It follows from the comparison principle that
$u_j\ge v_0,\ \forall j\in \N.$ Therefore by the
 general properties of psh functions $(u_j)_j$ is
 relatively compact in  $L^1_{loc}-$topology.
Let $\tilde u \in \m N^a (f) $ be any closter point of the
 sequence $u_j.$ Assume that $u_j \to \tilde u$
pointwise $d\lambda -$almost everywhere, here $d\lambda$ denotes the
Lebesgue measure. By Lemma 2.1 in
\cite{Ce1}, after extracting a subsequence if necessary,
 we have $u_j \to \tilde u $ $d\mu-$almost everywhere.
Then
$$
\tilde u = (\limsup_{j\to +\infty}u_j )^* = \lim_{j\to +\infty} (\sup_{k\ge j}u_k)^*.
$$
Now, consider the following auxiliary functions
$$
\tilde u_j= (\sup_{k\ge j}u_k)^* =
(\lim_{l\to +\infty} \sup_{l\ge k \ge j}u_k)^* =( \lim _{l\to +\infty}
\tilde{u}_j^l)^*.
$$
Observe that
$$
(dd^c \max(u_j , u_k))^n \ge \min (h_j , h_k)d\mu.
$$
Therefore
$$
(dd^c \tilde u_j )^n = \lim _{l\to +\infty}
(dd^c \tilde u_j^l )^n  \ge \lim_{l \to +\infty} \min_{l\ge k\ge j} h_k d\mu.
$$
We Let $j$ converges to $+\infty$ to  get
$$
(dd^c \tilde u)^n \ge h d\mu.
$$
Now, for the reverse inequality, pick $\f \in \m E_0$ a
 negative psh  function. For any $j\ge 1$ and
  since
$u_j \le \tilde u_j,$ we have by  integration by parts,
which is valid in $\m N^a(f)$ (cf \cite{ACCP}), that
$$
\int_\W -\f (dd^c u_j )^n  \ge \int_\W -\f (dd^c \tilde u_j )^n.
$$
Therefore
$$
\lim_{j\to +\infty}\int_\W \f h_j d\mu
 \le \lim_{j\to +\infty}\int_\W \f (dd^c \tilde u_j )^n
= \int_\W \f (dd^c \tilde u)^n.
$$
Together with the first inequality, we get
$$
 \int_\W \f (dd^c \tilde u)^n =  \int_\W \f h d\mu , \quad \forall \f \in \m E_0.
$$
The set of test functions $\m D(\W)\subset  \m E_0 - \m E_0 $
(cf Lemma 2.1 in \cite{Ce4}) therefore the equality holds for any
 $\f \in \m D(\W).$ Thence
$$
(dd^c \tilde u)^n = hd\mu = (dd^c u)^n .
$$
Uniqueness in the class $\m N^a(f)$ implies that $\tilde u = u$  which concludes the proof.
\end{proof}
 \begin{proof}[Proof of Main Theorem ]
Assume first that $F(t, .) \in L^1(d\mu).$ Then  $F(f, .) \in L^1(d\mu).$
 It follows from \cite{Ce8} and  \cite{ACCP} that the nonnegative measure $F(f, .) d \mu $
 is the Monge-Amp\`ere measure of a function
$u_0 $
from the class $\m F^a(f) .$ Then
$$
(dd^c u_0 )^n = F(f, .) d\mu \ge F(u_0 , .) d \mu .
$$
We denote by $\m A $ the set of all $u\in \m F^a(f) $ such that
 $u\ge u_0.$ The set $\m A$
is convex and compact with respect to $L^1(d\lambda)$-topology,
where $d\lambda $ denotes the Lebesgue measure
in $\C^n.$ Once more, by \cite{Ce8} (see also \cite{ACCP}), we have
for each $u\in \m A,$ there exists a unique
$\hat{u}  \in \m F^a(f)$ such that
$$
(dd^c \hat{u} )^n = F(u, .) d\mu.
$$
Since $\hat{u} \le f$ and $F$ is nondecreasing in the first variable then
$$
(dd^c \hat{u} )^n =  F(u, \cdot) d\mu \le F(f, \cdot) d\mu  =  (dd^c u_0)^n.
$$
The comparison principle yields that $\hat u \ge u\ge u_0, $
 hence  $\hat u \in \m A.$

We define the map $T : \m A \to \m A $ by $u \mapsto \hat{u}.$
By Schauder's fixed point theorem, we are done as
soon as we show that the map $T$ is continuous. Let $u_j \in \m A$ be a sequence
which converges towards $u \in \m A.$ By Lemma \ref{stability}, it's enough to show that
$F(u_j , .) d\mu \to F(u, .)d \mu .$ After extracting a subsequence, we may assume that
$u_j \to u $ $d\lambda $-a.e. Applying Lemma 2.1. in \cite{Ce1}, we get $u_j \to u $
$d\mu$-a.e. By Lebesgue  convergence theorem we have
$F(u_j , .) d\mu \to F(u, .)d \mu .$

We now proceed to complete the proof of the general case. Let us consider the set
$$
\m K : = \left\{  \f \in \m N^a(f);  (dd^c \f)^n \ge F(\f, .) d\mu\right\}.
$$
Claim 1. $\m K  $ is not empty: It follows from the monotonicity of $F$
$$
(dd^c v_0 +f)^n \ge (dd^c v_0)^n \ge F(v_0 , ) d\mu \ge F(v_0 + f , .) d\mu .
$$
 Then the function  $\f_0 :=v_0 + f$ belongs to $\m K .$

Let denote
$$
\m K_0 := \left\{ \f \in \m K ; \ \f \ge \f_0\right \}.
$$
Claim 2. $\m K_0 $  is stable by taking the maximum: Let $\f _1 , \f _2 \in \m K_0.$
It's clear that $\max(u_1, u_2) \ge \f_0.$
Since $\m N^a(f) $ is stable by taking the maximum then
$\max(u_1, u_2) \in \m N^a(f) .$ On the other hand, from
\cite{Dem}, we have
\begin{eqnarray*}
(dd^c \max (u_1, u_2))^n &\ge & \textbf{1}_{(u_1\ge u_2)}(dd^c u_1)^n +
\textbf{1}_{(u_1 < u_2)}(dd^c u_2)^n \\
&\ge & \textbf{1}_{(u_1\ge u_2)} F(u_1 , .) d\mu + \textbf{1}_{(u_1< u_2)}
F(u_2 , . ) d\mu \\
&\ge & F(\max(u_1, u_2), . ) d\mu.
\end{eqnarray*}
 Which implies that $\max(u_1, u_2) \in \m K_0.$

 Claim 3. $ \m K_0$ is compact in $L^1_{loc} (\W):$  It's enough to prove that it's closed.
 Let $\f _j \in \m K_0$ be a sequence converging towards
 $ \f \in \m N^a(f).$ The limit function
 is given by $\f = (\limsup_{j\to \infty} \f_j)^*.$ Then
 $\f _0 \le \f \le f.$ The continuity of the  complex Monge-
Amp\`ere operator and the properties of $F$ yield
\begin{eqnarray*}
(dd^c \f )^n &= &\lim_{j\to +\infty} (dd^c  \sup _{k\ge j} \f_k )^n\\
             &= & \lim_{j\to +\infty} \lim_{l\to +\infty}
             (dd^c  \max_{l\ge k\ge j} \f_k )^n\\
             &\ge & \lim_{j\to +\infty} \lim_{l\to +\infty} F(\max _{l\ge k \ge j}
             \f_k , .) d\mu .
\end{eqnarray*}
Therefore $\f \in \m K_0.$

Consider the following upper envelope
$$
\phi(z) : = \sup \left\lbrace \f (z) ; \f \in \m K_0 \right\rbrace.
$$
Notice that in order to get a psh function $\phi$ we should a priori replace $\phi$
by its upper semi-continuous regularization
$\phi^*(z) := \limsup_{\zeta \to z} \phi (\zeta )$ but since $\phi ^* \in \m K_0$ as well
$ \phi^*$ contributes to the envelope (i.e. $\phi ^* \in \m K _0$)  and then $\phi = \phi^*.$

Claim 4. $\phi $ is solution to Monge-Amp\`ere equation (\ref{local}):
It follows from Choquet's Lemma that there exists a sequence
$\phi_j \in \m K_0 $ such that
$$
\phi =(\limsup_{j\to +\infty} \phi _ j)^*.
$$
Since the class $\m K_0$ is stable under taking the maximum, we can assume that $\phi_j$
is nondecreasing.
 We use the classical balayage procedure to prove that $\phi$
 is actually a solution of (\ref{local}). Pick $\textbf{B}\Sub \W$
  a ball and define the function
$$
\phi _j^B (z) := \sup \left\lbrace  v(z) ; \ v^* \le \phi_j \ \text{on\ }\
 \partial \textbf{B}, \ v\in PSH(\textbf{B})\right\rbrace, \ \forall  z\in \textbf{B}.
$$
By the first case,  there exists a function
$\tilde{\phi}_j \in \m F^a (\phi_j^B , \textbf{B}) $ solution to the following equation
$$
(dd^c \tilde{\phi}_j ) ^n = \textbf{1}_{\textbf{B}} F(\tilde{\phi}_j , .) d\mu.
$$
In fact, the function $\tilde \phi_j$ is given by the following upper envelope
$$
\tilde \phi_j = \sup\left\lbrace  w\in \m E(\textbf{B}) ; \ w \le \phi _j^B  \ \text{and\  }
\ (dd^c w )^n \ge F(w, .) d\mu \right\rbrace.
$$
Indeed, if we denote by $g$ the right hand side function, then
 $\tilde{\phi}_j \le g \le \phi_j^B.$ Hence $g\in \m F^a(\phi_j^B, \textbf{B}).$
It follows by Lemma 3.3 in \cite{ACCP} that
\begin{equation}\label{1}
\int_\W  \chi (dd^c \tilde{\phi}_j)^n \le
\int _\W \chi (dd^c g)^n \ , \quad \forall \chi \in \m E_0.
\end{equation}
 On the  other hand, as before, we have
$g = (\lim g_k)^*$ where  $g_k \in \m E(\textbf{B})$ is a nondecreasing sequence satisfying
$  \phi_j^B \ge g_k \ge \phi_j $  and $(dd^c g_k )^n \ge F(g_k , \cdot) d \mu .$ Therefore
 $(dd^c g )^n \ge F(g , \cdot) d \mu .$
Then
\begin{equation}\label{2}
(dd^c \tilde{\phi}_j )^n = F( \tilde{\phi}_j, .) d\mu \le
F(g, . ) d\mu  \le (dd^c g)^n.
\end{equation}
Combining (\ref{1}) and (\ref{2}), one get
$$
(dd^c \tilde{\phi}_j )^n = (dd^c g)^n,
$$
therefore, by the comparison principle, we have $\tilde{\phi}_j =g.$

Now, for $j\in \N,$ let consider the function $\psi_j $ defined on $\W$ by
$$
\psi_j (z)  =
\left\lbrace
\begin{array}{cc}
 \tilde{\phi}_j (z) \ \text{if} \  z\in \textbf{B}\\
 \phi_j (z)  \  \text{if}  \  z \not \in \textbf{B}
\end{array}
\right.
.
$$
On $\textbf{B}$ we have $ \phi_j \le \tilde{\phi }_j \le  \phi _j^B \le f$
and on $\W \setminus \textbf{B} $ we have $\tilde{\phi }_j  = \phi _ j \le f .$ Hence
$\psi_j \in \m N^a(f).$ From the definition of $\psi_j $, we deduce that
 $ (dd^c \psi_j )^n \ge F( \psi_j , \cdot) d\mu.$
Therefore
 $\psi_j \in \m K_0$  and
$$
\phi =(\lim_{j\to +\infty} \psi_j)^*.
$$
 Since the complex Monge-Amp\`ere operator
 is continuous under monotonic sequences
and $\textbf{B}$ is arbitrarily chosen, to conclude the proof of the claim  it's enough to observe that
 the sequence
$\psi_j$ is nondecreasing.

Uniqueness follows in a classical way from the comparison principle and
the monotonicity of the function $F$.
 Indeed, assume that there exist two solutions
 $\f _1$ and $\f_2$ in $\m N^a(f)$ such that
 $$
 (dd^c \f_i)^n = F(\f_i,.) d \mu ,\quad i=1, 2.
 $$
 Since the function $F$ is nondecreasing in the first variable, then
$$
 F(\f_1 , .)d\mu \le  F(\f_2 , .)d\mu \quad \text{on}\quad (\f_1<\f_2).
$$
On the other hand, by the comparison principle we have
 \begin{multline*}
 \int_{(\f_1<\f_2)} F(\f_2 , .)d\mu =\int_{(\f_1<\f_2)}  (dd^c \f_2)^n\le \\
 \int_{(\f_1<\f_2)} (dd^c \f_1)^n=
  \int_{(\f_1<\f_2)} F(\f_2 , .)d\mu.
 \end{multline*}
  Therefore
  $$
 F(\f_1 , .)d\mu = F(\f_2 , .)d\mu \quad \text{on}\ (\f_1 < \f_2) .
  $$
 In the same way, we get the equality on $(\f_1>\f_2)$ and then on $\W.$ Hence
 $(dd^c \f_1 )^n =(dd^c \f_2)^n \ \text{on}\ \W.$ Therefore  uniqueness in the class $\m N^a(f)$ yields
  $\f_1 = \f_2$  and the proof is completed.
\end{proof}
\begin{rqes}
1-- We have no precise knowledge  when the subsolution of the equation
(\ref{local} ) exists. However,  if there
exists a negative function $\psi \in PSH(\W)$ such that
$$
\int _\W -\psi F(0, .) d\mu < \infty ,
$$
then (\ref{local} ) admits a subsolution $v\in \m N^a.$ This is an immediate consequence of
Proposition 5.2 in \cite{Ce8}.

2-- The condition 2 in the Theorem  is necessary.

\end{rqes}


\subsection*{Acknowledgements}
The author is grateful to the referee for his comments and suggestions.
This note was written during the
author's visit to Institut de Math\'ematiques de Toulouse .
He  likes to thank Vincent Guedj and Ahmed Zeriahi  for fruitful discussions
and the warm hospitality.

\end{document}